\documentclass[12pt]{article}
\usepackage{amsmath, mathrsfs}
\usepackage{amssymb,amsmath,amsfonts}
\usepackage{hyperref}
\usepackage{color}
\newtheorem{theorem}{Theorem}[section]
\newtheorem{lem}[theorem]{Lemma}
\newtheorem{prop}[theorem]{Proposition}
\newtheorem{cor}[theorem]{Corollary}
\newtheorem{defn}[theorem]{Definition}
\newenvironment{defn-new}{\begin{defn} \em}{\end{defn}}
\newtheorem{rem}[theorem]{Remark}
\newenvironment{rem-new}{\begin{rem} \em}{\end{rem}}
\newtheorem{ex}[theorem]{Example}
\newenvironment{ex-new}{\begin{ex} \em}{\end{ex}}
\newtheorem{exer}[theorem]{Exercise}
\newenvironment{exer-new}{\begin{exer} \em}{\end{exer}}
\newtheorem{agr}[theorem]{Agreement}
\newenvironment{agr-new}{\begin{agr} \em}{\end{agr}}
\newtheorem{pbm}[theorem]{Problem}
\newenvironment{pbm-new}{\begin{pbm} \em}{\end{pbm}}

\makeatletter \@addtoreset{equation}{section} \makeatother

\begin{document}

\begin{center}
{\Large {\bf Indefinite almost paracontact metric manifolds}}\bigskip
\bigskip

Mukut Mani Tripathi, Erol K\i l\i \c{c}, Selcen Y\"{u}ksel
Perkta\c{s} and Sad\i k Kele\c{s}\bigskip \bigskip
\end{center}

\noindent {\bf Abstract.} In this paper we introduce the concept of $\left(
\varepsilon \right) $-almost paracontact manifolds, and in particular, of $%
\left( \varepsilon \right) $-para Sasakian manifolds. Several examples are
presented. Some typical identities for curvature tensor and Ricci tensor of $%
\left( \varepsilon \right) $-para Sasakian manifolds are obtained. We prove
that if a semi-Riemannian manifold is one of flat, proper recurrent or
proper Ricci-recurrent, then it can not admit an $\left( \varepsilon \right)
$-para Sasakian structure. We show that, for an $\left( \varepsilon \right) $%
-para Sasakian manifold, the conditions of being symmetric, semi-symmetric
or of constant sectional curvature are all identical. It is shown that a
symmetric spacelike (resp. timelike) $\left( \varepsilon \right) $-para
Sasakian manifold $M^{n}$ is locally isometric to a pseudohyperbolic space $%
H_{\nu }^{n}\left( 1\right) $ (resp. pseudosphere $S_{\nu }^{n}\left(
1\right) $). In last, it is proved that for an $\left( \varepsilon \right) $%
-para Sasakian manifold, the conditions of being Ricci-semisymmetric,
Ricci-symmetric and Einstein are all identical. \medskip

 \medskip


\section{Introduction\label{sect-intro}}

In 1976, an almost paracontact structure $(\varphi ,\xi ,\eta )$ satisfying $%
\varphi ^{2}=I-\eta \otimes \xi $ and $\eta (\xi )=1$ on a differentiable
manifold, was introduced by S\={a}to \cite{Sato-76}. The structure is an
analogue of the almost contact structure \cite%
{Sasaki-60-Tohoku,Blair-02-book} and is closely related to almost product
structure (in contrast to almost contact structure, which is related to
almost complex structure). An almost contact manifold is always
odd-dimensional but an almost paracontact manifold could be even-dimensional
as well. In 1969, T. Takahashi \cite{Takahashi-69-Tohoku-1} introduced
almost contact manifolds equipped with associated pseudo-Riemannian metrics.
In particular, he studied Sasakian manifolds equipped with an associated
pseudo-Riemannian metric. These indefinite almost contact metric manifolds
and indefinite Sasakian manifolds are also known as $\left( \varepsilon
\right) $-almost contact metric manifolds and $\left( \varepsilon \right) $%
-Sasakian manifolds respectively \cite%
{Bej-Dug-93,Duggal-90-IJMMS,Duggal-Sahin-07}. Also, in 1989, K. Matsumoto
\cite{Mat-89} replaced the structure vector field $\xi $ by $-\,\xi $ in an
almost paracontact manifold and associated a Lorentzian metric with the
resulting structure and called it a Lorentzian almost paracontact manifold.
\medskip

An $(\varepsilon)$-Sasakian manifold is always odd-dimensional. Recently, we
have observed that there does not exist a lightlike surface in a $3$%
-dimensional $(\varepsilon)$-Sasakian manifold. On the other hand,
in a Lorentzian almost paracontact manifold given by Matsumoto, the
semi-Riemannian metric has only index $1$ and the structure vector
field $\xi $ is always timelike. These circumstances motivate us to
associate a semi-Riemannian metric, not necessarily Lorentzian, with
an almost paracontact structure, and we shall call this indefinite
almost paracontact metric structure an $\left( \varepsilon
\right) $-almost paracontact structure, where the structure vector field $%
\xi $ will be spacelike or timelike according as $\varepsilon =1$ or $%
\varepsilon =-1$. \medskip

In this paper we initiate study of $\left( \varepsilon \right) $-almost
paracontact manifolds, and in particular, $\left( \varepsilon \right) $-para
Sasakian manifolds. The paper is organized as follows. Section~\ref%
{sect-eps-ACMM} contains basic definitions and some examples of $\left(
\varepsilon \right) $-almost paracontact manifolds. In section~\ref%
{sect-apcm-normal}, some properties of normal almost paracontact structures
are discussed. Section~\ref{sect-eps-spc} contains definitions of an $\left(
\varepsilon \right) $-paracontact structure and an $\left( \varepsilon
\right) $-$s$-paracontact structure. A typical example of an $\left(
\varepsilon \right) $-$s$-paracontact structure is also presented. In
section~\ref{sect-eps-PS}, we introduce the notion of an $\left( \varepsilon
\right) $-para Sasakian structure and study some of its basic properties. We
find some typical identities for curvature tensor and Ricci tensor. We prove
that if a semi-Riemannian manifold is one of flat, proper recurrent or
proper Ricci-recurrent, then it can not admit an $\left( \varepsilon \right)
$-para Sasakian structure. We show that, for an $\left( \varepsilon \right) $%
-para Sasakian manifold, the conditions of being symmetric, semi-symmetric
or of constant sectional curvature are all identical. More specifically, it
is shown that a symmetric spacelike $\left( \varepsilon \right) $-para
Sasakian manifold $M^{n}$ is locally isometric to a pseudohyperbolic space $%
H_{\nu }^{n}\left( 1\right) $ and a symmetric timelike $\left( \varepsilon
\right) $-para Sasakian manifold $M^{n}$ is locally isometric to a
pseudosphere $S_{\nu }^{n}\left( 1\right) $. In last, it is proved that for
an $\left( \varepsilon \right) $-para Sasakian manifold, the conditions of
being Ricci-semisymmetric, Ricci-symmetric and Einstein are all identical.
Unlike $3$-dimensional $(\varepsilon )$-Sasakian manifold, which cannot
possess a lightlike surface, the study of lightlike surfaces of $3$%
-dimensional $\left( \varepsilon \right) $-para Sasakian manifolds will be
presented in a forthcoming paper.

\section{$(\protect\varepsilon )$-almost paracontact metric manifolds\label%
{sect-eps-ACMM}}

Let $M$ be an almost paracontact manifold \cite{Sato-76} equipped with an
almost paracontact structure $(\varphi ,\xi ,\eta )$ consisting of a tensor
field $\varphi $ of type $(1,1)$, a vector field $\xi $ and a $1$-form $\eta
$ satisfying
\begin{equation}
\varphi ^{2}=I-\eta \otimes \xi ,  \label{eq-phi-eta-xi}
\end{equation}%
\begin{equation}
\eta (\xi )=1,  \label{eq-eta-xi}
\end{equation}%
\begin{equation}
\varphi \xi =0,  \label{eq-phi-xi}
\end{equation}%
\begin{equation}
\eta \circ \varphi =0.  \label{eq-eta-phi}
\end{equation}%
It is easy to show that the relation (\ref{eq-phi-eta-xi}) and one of the
three relations (\ref{eq-eta-xi}), (\ref{eq-phi-xi}) and (\ref{eq-eta-phi})
imply the remaining two relations of (\ref{eq-eta-xi}), (\ref{eq-phi-xi})
and (\ref{eq-eta-phi}). On an $n$-dimensional almost paracontact manifold,
one can easily obtain
\begin{equation}
\varphi ^{3}-\varphi =0,  \label{eq-phi3+phi}
\end{equation}%
\begin{equation}
{\rm rank}\left( \varphi \right) =n-1.  \label{eq-rank-phi}
\end{equation}%
The equation $(\ref{eq-phi3+phi})$ gives an $f(3\,,-1)${\em -structure} \cite%
{Singh-Vohra-72}. \medskip

Throughout the paper, by a semi-Riemannian metric \cite{ONeill-83} on a
manifold $M$, we understand a non-degenerate symmetric tensor field $g$ of
type $\left( 0,2\right) $. In particular, if its index is $1$, it becomes a
Lorentzian metric \cite{Beem-Ehrlich-81}. A sufficient condition for the
existence of a Riemannian metric on a differentiable manifold is
paracompactness. The existence of Lorentzian or other semi-Riemannian
metrics depends upon other topological properties. For example, on a
differentiable manifold, the following statements are equivalent: (1) there
exits a Lorentzian metric on $M$, (2) there exists a non vanishing vector
field on $M$, (3) either $M$ is non compact, or $M$ is compact and has Euler
number $\chi (M)=0$. Also for instance, the only compact surfaces that can
be made Lorentzian surfaces are the tori and Klein bottles, and a sphere $%
S^{n}$ admits a Lorentzian metric if and only if $n$ is odd $\geq 3$.
\medskip

Now, we give the following:

\begin{defn-new}
\label{defn-ind-apcm} Let $M$ be a manifold equipped with an almost
paracontact structure $(\varphi ,\xi ,\eta )$. Let $g$ be a semi-Riemannian
metric with ${\rm index}(g)=\nu $ such that
\begin{equation}
g\left( \varphi X,\varphi Y\right) =g\left( X,Y\right) -\varepsilon \eta
(X)\eta \left( Y\right) ,\qquad X,Y\in TM,  \label{eq-metric-1}
\end{equation}%
where $\varepsilon =\pm 1$. Then we say that $M$ is an $\left( \varepsilon
\right) $-{\em almost paracontact metric manifold} equipped with an $\left(
\varepsilon \right) ${\em -almost paracontact metric structure} $(\varphi
,\xi ,\eta ,g,\varepsilon )$. In particular, if ${\rm index}(g)=1$, then an $%
(\varepsilon )$-almost paracontact metric manifold will be called a {\em %
Lorentzian almost paracontact manifold}. In particular, if the metric $g$ is
positive definite, then an $(\varepsilon )$-almost paracontact metric
manifold is the usual {\em almost paracontact metric manifold} \cite{Sato-76}%
.
\end{defn-new}

The equation (\ref{eq-metric-1}) is equivalent to
\begin{equation}
g\left( X,\varphi Y\right) =g\left( \varphi X,Y\right)  \label{eq-metric-2}
\end{equation}%
along with
\begin{equation}
g\left( X,\xi \right) =\varepsilon \eta (X)  \label{eq-metric-3}
\end{equation}%
for all $X,Y\in TM$. From (\ref{eq-metric-3}) it follows that
\begin{equation}
g\left( \xi ,\xi \right) =\varepsilon ,  \label{eq-g(xi,xi)}
\end{equation}%
that is, the structure vector field $\xi $ is never lightlike. Since $g$ is
non-degenerate metric on $M$ and $\xi $ is non-null, therefore the
paracontact distribution
\[
D=\{X\in TM:\eta \left( X\right) =0\}
\]%
is non-degenerate on $M$. \medskip

\begin{defn-new}
Let $(M,\varphi ,\xi ,\eta ,g,\varepsilon )$ be an $(\varepsilon )$-almost
paracontact metric manifold (resp. a Lorentzian almost paracontact
manifold). If $\varepsilon =1$, then $M$ will be said to be a {\em spacelike}
$(\varepsilon )$-{\em almost paracontact metric manifold} (resp. a {\em %
spacelike Lorentzian almost paracontact manifold}). Similarly, if $%
\varepsilon =-\,1$, then $M$ will be said to be a {\em timelike} $%
(\varepsilon )$-{\em almost paracontact metric manifold} (resp. a {\em %
timelike Lorentzian almost paracontact manifold}).
\end{defn-new}

Note that a timelike Lorentzian almost paracontact structure is a Lorentzian
almost paracontact structure in the sense of Mihai and Rosca \cite%
{Mih-Rosca-92,Mat-Mih-Rosca-95}, which differs in the sign of the structure
vector field of the Lorentzian almost paracontact structure given by
Matsumoto \cite{Mat-89}.

\begin{ex-new}
\label{ex-eps-apcm-1} Let ${\Bbb R}^{3}$\ be the $3$-dimensional real number
space with a coordinate system $\left( x,y,z\right) $. We define
\[
\eta =dy\ ,\qquad \xi =\frac{\partial }{\partial y}\,,
\]%
\[
\varphi \left( \frac{\partial }{\partial x}\right) =\frac{\partial }{%
\partial z}\ ,\qquad \varphi \left( \frac{\partial }{\partial y}\right) =0\
,\qquad \varphi \left( \frac{\partial }{\partial z}\right) =\frac{\partial }{%
\partial x}\,,
\]%
\[
g_{1}=\left( dx\right) ^{2}-\left( dy\right) ^{2}+\left( dz\right) ^{2}\,,
\]%
\[
g_{2}=-\left( dx\right) ^{2}+\left( dy\right) ^{2}-\left( dz\right) ^{2}\,.
\]%
Then the set $\left( \varphi ,\xi ,\eta ,g_{1}\right) $ is a timelike
Lorentzian almost paracontact structure, while the set $\left( \varphi ,\xi
,\eta ,g_{2}\right) $ is a spacelike $\left( \varepsilon \right) $-almost
paracontact metric structure. We note that {\rm index}$\left( g_{1}\right)
=1 $ and {\rm index}$\left( g_{2}\right) =2$.
\end{ex-new}

\begin{ex-new}
\label{ex-eps-apcm-2} Let ${\Bbb R}^{3}$\ be the $3$-dimensional real number
space with a coordinate system $\left( x,y,z\right) $. We define
\[
\eta =dz-y\,dx\ ,\qquad \xi =\frac{\partial }{\partial z}\,,
\]%
\[
\varphi \left( \frac{\partial }{\partial x}\right) =-\,\frac{\partial }{%
\partial x}-y\frac{\partial }{\partial z}\ ,\qquad \varphi \left( \frac{%
\partial }{\partial y}\right) =-\,\frac{\partial }{\partial y}\,,\qquad
\varphi \left( \frac{\partial }{\partial z}\right) =0\ ,
\]%
\[
g_{1}=\left( dx\right) ^{2}+\left( dy\right) ^{2}-\eta \otimes \eta \,,
\]%
\[
g_{2}=\left( dx\right) ^{2}+\left( dy\right) ^{2}+\left( dz\right)
^{2}-y\left( dx\otimes dz+dz\otimes dx\right) \,,
\]%
\[
g_{3}=-\left( dx\right) ^{2}+\left( dy\right) ^{2}+\left( dz\right)
^{2}-y\left( dx\otimes dz+dz\otimes dx\right) \,.
\]%
Then, the set $\left( \varphi ,\xi ,\eta \right) $ is an almost paracontact
structure in ${\Bbb R}^{3}$. The set $\left( \varphi ,\xi ,\eta
,g_{1}\right) $ is a timelike Lorentzian almost paracontact structure.
Moreover, the trajectories of the timelike structure vector $\xi $ are
geodesics. The set $\left( \varphi ,\xi ,\eta ,g_{2}\right) $ is a spacelike
Lorentzian almost paracontact structure. The set $\left( \varphi ,\xi ,\eta
,g_{3}\right) $ is a spacelike $(\varepsilon )$-almost paracontact metric
structure $\left( \varphi ,\xi ,\eta ,g_{3},\varepsilon \right) $ with {\rm %
index}$\left( g_{3}\right) =2$.
\end{ex-new}

\begin{ex-new}
\label{ex-eps-apcm-3} Let ${\Bbb R}^{5}$\ be the $5$-dimensional real number
space with a coordinate system $\left( x,y,z,t,s\right) $. Defining
\[
\eta =ds-ydx-tdz\ ,\qquad \xi =\frac{\partial }{\partial s}\ ,
\]%
\[
\varphi \left( \frac{\partial }{\partial x}\right) =-\,\frac{\partial }{%
\partial x}-y\frac{\partial }{\partial s}\ ,\qquad \varphi \left( \frac{%
\partial }{\partial y}\right) =-\,\frac{\partial }{\partial y}\ ,
\]%
\[
\varphi \left( \frac{\partial }{\partial z}\right) =-\,\frac{\partial }{%
\partial z}-t\frac{\partial }{\partial s}\ ,\qquad \varphi \left( \frac{%
\partial }{\partial t}\right) =-\,\frac{\partial }{\partial t}\ ,\qquad
\varphi \left( \frac{\partial }{\partial s}\right) =0\ ,
\]%
\[
g_{1}=\left( dx\right) ^{2}+\left( dy\right) ^{2}+\left( dz\right)
^{2}+\left( dt\right) ^{2}-\eta \otimes \eta \ ,
\]%
\begin{eqnarray*}
g_{2} &=&-\,\left( dx\right) ^{2}-\left( dy\right) ^{2}+\left( dz\right)
^{2}+\left( dt\right) ^{2}+\left( ds\right) ^{2} \\
&&-\,t\left( dz\otimes ds+ds\otimes dz\right) -y\left( dx\otimes
ds+ds\otimes dx\right) ,
\end{eqnarray*}%
the set $(\varphi ,\xi ,\eta ,g_{1})$\ becomes a timelike Lorentzian almost
paracontact structure in ${\Bbb R}^{5}$, while the set $(\varphi ,\xi ,\eta
,g_{2})$\ is a spacelike $\left( \varepsilon \right) $-almost paracontact
structure. Note that {\rm index}$\left( g_{2}\right) =3$.
\end{ex-new}

The Nijenhuis tensor $\left[ J,J\right] $ of a tensor field $J$ of type $%
\left( 1,1\right) $ on a manifold $M$ is a tensor field of type $\left(
1,2\right) $ defined by
\begin{equation}
\left[ J,J\right] (X,Y)\equiv J^{2}\left[ X,Y\right] +\left[ JX,JY\right] -J%
\left[ JX,Y\right] -J\left[ X,JY\right]  \label{eq-Nij-tens}
\end{equation}%
for all $X,Y\in TM$. If $M$ admits a tensor field $J$ of type $\left(
1,1\right) $ satisfying
\begin{equation}
J^{2}=I,  \label{eq-ap-str}
\end{equation}%
then it is said to be an {\em almost product manifold} equipped with an {\em %
almost product structure} $J$. An almost product structure is {\em integrable%
} if its Nijenhuis tensor vanishes. For more details we refer to \cite%
{Yano-65-book}.

\begin{ex-new}
Let $\left( M^{n},J,G\right) $ be a semi-Riemannian almost product manifold,
such that
\[
J^{2}=I,\qquad G\left( JX,JY\right) =G\left( X,Y\right) .
\]%
Consider the product manifold $M^{n}\times {\Bbb R}$. A vector field on $%
M^{n}\times {\Bbb R}$ can be represented by $\left( X,f\frac{d}{dt}\right) $%
, where $X$ is tangent to $M$, $f$ a smooth function on $M^{n}\times {\Bbb R}
$ and $t$ the coordinates of ${\Bbb R}$. On $M^{n}\times {\Bbb R}$ we define
\[
\eta =dt,\qquad \xi =\frac{d}{dt},\qquad \varphi \left( \left( X,f\frac{d}{dt%
}\right) \right) =JX,
\]%
\[
g\left( \left( X,f\frac{d}{dt}\right) ,\left( Y,h\frac{d}{dt}\right) \right)
=G\left( X,Y\right) +\varepsilon fh.
\]%
Then $\left( \varphi ,\xi ,\eta ,g,\varepsilon \right) $ is an $\left(
\varepsilon \right) $-almost paracontact metric structure on the product
manifold $M^{n}\times {\Bbb R}$.
\end{ex-new}

\begin{ex-new}
Let $\left( M,\psi ,\xi ,\eta ,g,\varepsilon \right) $ be an $\left(
\varepsilon \right) $-almost contact metric manifold. If we put $\varphi
=\psi ^{2}$, then $\left( M,\varphi ,\xi ,\eta ,g,\varepsilon \right) $ is
an $\left( \varepsilon \right) $-almost paracontact metric manifold.
\end{ex-new}

\section{Normal almost paracontact manifolds\label{sect-apcm-normal}}

Let $M$ be an almost paracontact manifold with almost paracontact structure $%
\left( \varphi ,\xi ,\eta \right) $ and consider the product manifold $%
M\times {\Bbb R}$, where ${\Bbb R}$ is the real line. A vector field on $%
M\times {\Bbb R}$ can be represented by $\left( X,f\frac{d}{dt}\right) $,
where $X$ is tangent to $M$, $f$ a smooth function on $M\times {\Bbb R}$ and
$t$ the coordinates of ${\Bbb R}$. For any two vector fields $\left( X,f%
\frac{d}{dt}\right) $ and $\left( Y,h\frac{d}{dt}\right) $, it is easy to
verify the following
\begin{equation}
\left[ \left( X,f\frac{d}{dt}\right) ,\left( Y,h\frac{d}{dt}\right) \right]
=\left( \left[ X,Y\right] ,\left( Xh-Yf\right) \frac{d}{dt}\right) .
\label{eq-Lie-bracket-prod-mfd}
\end{equation}

\begin{defn-new}
If the induced almost product structure $J$ on $M\times {\Bbb R}$ defined by
\begin{equation}
J\left( X,f\frac{d}{dt}\right) \equiv \left( \varphi X+f\xi \,,\,\eta \left(
X\right) \frac{d}{dt}\right)  \label{eq-induced-J}
\end{equation}%
is integrable, then we say that the almost paracontact structure $\left(
\varphi ,\xi ,\eta \right) $ is {\em normal}.
\end{defn-new}

This definition is conformable with the definition of normality given in
\cite{BHM-81}. As the vanishing of the Nijenhuis tensor $\left[ J,J\right] $
is a necessary and sufficient condition for the integrability of the almost
product structure $J$, we seek to express the conditions of normality in
terms of the Nijenhuis tensor $\left[ \varphi ,\varphi \right] $ of $\varphi
$. In view of (\ref{eq-Nij-tens}), (\ref{eq-induced-J}), (\ref%
{eq-Lie-bracket-prod-mfd}) and (\ref{eq-phi-eta-xi})--(\ref{eq-eta-phi}) we
have
\begin{eqnarray*}
&&\left[ J,J\right] \left( \left( X,f\frac{d}{dt}\right) ,\left( Y,h\frac{d}{%
dt}\right) \right) \medskip \\
&&\quad =\left( \!\!\frac{{}}{{}}\left[ \varphi ,\varphi \right]
(X,Y)-2d\eta (X,Y)\xi \right. -\ h\left( {\pounds }_{\xi }\varphi \right)
X+f\left( {\pounds }_{\xi }\varphi \right) Y,\medskip \\
&&\qquad \left. \left\{ \left( {\pounds }_{\varphi X}\eta \right) Y-\left( {%
\pounds }_{\varphi Y}\eta \right) X-h\left( {\pounds }_{\xi }\eta \right)
X+f\left( {\pounds }_{\xi }\eta \right) Y\right\} \frac{d}{dt}\right)
,\smallskip
\end{eqnarray*}%
where ${\pounds }_{X}$ denotes the Lie derivative with respect to $X$. Since
$\left[ J,J\right] $ is skew symmetric tensor field of type $\left(
1,2\right) $, it suffices to compute $\left[ J,J\right] \left( \left(
X,0\right) ,\left( Y,0\right) \right) $ and $\left[ J,J\right] \left( \left(
X,0\right) ,\left( 0,\frac{d}{dt}\right) \right) $. Thus we have
\begin{eqnarray*}
\left[ J,J\right] \left( \left( X,0\right) ,\left( Y,0\right) \right)
&=&\left( \!\!\frac{{}}{{}}\left[ \varphi ,\varphi \right] (X,Y)-2d\eta
(X,Y)\xi \right. , \\
&&\left. \quad \left( \left( {\pounds }_{\varphi X}\eta \right) Y-\left( {%
\pounds }_{\varphi Y}\eta \right) X\right) \frac{d}{dt}\right) ,
\end{eqnarray*}%
\[
\left[ J,J\right] \left( \left( X,0\right) ,\left( 0,\frac{d}{dt}\right)
\right) =-\left( \left( {\pounds }_{\xi }\varphi \right) X\,,\,\left( \left(
{\pounds }_{\xi }\eta \right) X\right) \frac{d}{dt}\right) .
\]

We are thus led to define four types of tensors ${\overset{1}{N}}$, $\overset%
{2}{N}$, $\overset{3}{N}$ and $\overset{4}{{N}}$ respectively by (see also
\cite{Sato-76})
\begin{equation}
{\overset{1}{N}}\equiv \left[ \varphi ,\varphi \right] -2d\eta \otimes \xi ,
\label{eq-N-1}
\end{equation}%
\begin{equation}
{\overset{2}{N}}\equiv \left( {\pounds }_{\varphi X}\eta \right) Y-\left( {%
\pounds }_{\varphi Y}\eta \right) X,  \label{eq-N-2}
\end{equation}%
\begin{equation}
{\overset{3}{N}}\equiv {\pounds }_{\xi }\varphi ,  \label{eq-N-3}
\end{equation}%
\begin{equation}
{\overset{4}{N}}\equiv {\pounds }_{\xi }\eta .  \label{eq-N-4}
\end{equation}%
Thus the almost paracontact structure $\left( \varphi ,\xi ,\eta \right) $
will be normal if and only if the tensors defined by (\ref{eq-N-1})--(\ref%
{eq-N-4}) vanish identically. \medskip

Taking account of (\ref{eq-phi-eta-xi})--(\ref{eq-phi3+phi}) and (\ref%
{eq-N-1})--(\ref{eq-N-4}) it is easy to obtain the following:

\begin{lem}
\label{lem-eps-normal} Let $M$ be an almost paracontact manifold with an
almost paracontact structure $\left( \varphi ,\xi ,\eta \right) $. Then
\begin{equation}
{\overset{4}{N}}\left( X\right) =2d\eta \left( \xi ,X\right) ,
\label{eq-N4(X)}
\end{equation}%
\begin{equation}
{\overset{2}{N}}\left( X,Y\right) =2\left( d\eta \left( \varphi X,Y\right)
+d\eta \left( X,\varphi Y\right) \right) ,  \label{eq-N2(X,Y)}
\end{equation}%
\begin{equation}
{\overset{1}{N}}\left( X,\xi \right) =-\,{\overset{3}{N}}\left( \varphi
X\right) =-\,\left[ \xi ,X\right] +\varphi \left[ \xi ,\varphi X\right] +\xi
\left( \eta \left( X\right) \right) \xi ,  \label{eq-N1(X,xi)}
\end{equation}%
\begin{equation}
{\overset{1}{N}}\left( \varphi X,Y\right) =-~\varphi \left[ \varphi ,\varphi %
\right] \left( X,Y\right) -{\overset{2}{N}}\left( X,Y\right) \xi -\eta
\left( X\right) {\overset{3}{N}}\left( Y\right) .  \label{eq-N1(phiX,Y)}
\end{equation}%
Consequently,
\begin{equation}
\overset{2}{N}\left( X,\varphi Y\right) =2\left( d\eta \left( \varphi
X,\varphi Y\right) +d\eta \left( X,Y\right) \right) +\eta \left( Y\right) {%
\overset{4}{N}}\left( X\right) ,  \label{eq-N2(X,phiY)}
\end{equation}%
\begin{equation}
\overset{4}{{N}}\left( X\right) =\eta ({\overset{1}{N}}\left( X,\xi \right)
)=\overset{2}{N}\left( \xi ,\varphi X\right) =-\,\eta (\overset{3}{N}\left(
\varphi X\right) ),  \label{eq-N4(X)-1}
\end{equation}%
\begin{equation}
{\overset{4}{N}}\left( \varphi X\right) =-\,\eta \left( \left[ \xi ,\varphi X%
\right] \right) =-\,\eta (\overset{3}{N}\left( X\right) ),
\label{eq-N4(phiX)}
\end{equation}%
\begin{equation}
\varphi ({\overset{1}{N}}\left( X,\xi \right) )=\overset{3}{N}\left(
X\right) +\overset{4}{{N}}\left( \varphi X\right) \xi ,
\label{eq-phi(N1(X,xi))}
\end{equation}%
\begin{equation}
\eta ({\overset{1}{N}}\left( \varphi X,Y\right) )=-\,\overset{2}{N}\left(
X,Y\right) +\eta \left( X\right) \overset{4}{{N}}\left( \varphi Y\right) .
\label{eq-eta(N1(phiX,Y))}
\end{equation}
\end{lem}

From (\ref{eq-N4(X)-1}), it follows that if $\overset{2}{N}$ or
$\overset{3}{N}$ vanishes then $\overset{4}{N}$ vanishes. In view of
(\ref{eq-N4(X)-1}), (\ref{eq-phi(N1(X,xi))}) and
(\ref{eq-eta(N1(phiX,Y))}), we can state the following

\begin{theorem}
\label{th-nor-al-cont-str} If in an almost paracontact manifold $M$, ${%
\overset{1}{N}}$ vanishes then $\overset{2}{N}$, $\overset{3}{N}$ and $%
\overset{4}{{N}}$ vanish identically. Hence, the almost paracontact
structure is normal if and only if ${\overset{1}{N}}=0$.
\end{theorem}

Some equations given in Lemma~\ref{lem-eps-normal} are also in
\cite{Sato-76}. First part of the Theorem~\ref{th-nor-al-cont-str}
is given as Theorem~3.4 of \cite{Sato-76}. Now, we find a necessary
and sufficient condition for the vanishing of $\overset{2}{N}$ in
the following

\begin{prop}
The tensor $\overset{2}{N}$ vanishes if and only if
\begin{equation}
d\eta \left( \varphi X,\varphi Y\right) =-\,d\eta (X,Y).
\label{eq-d-eta(phiX,phiY)}
\end{equation}
\end{prop}

\noindent {\em Proof.} The necessary part follows from (\ref{eq-N2(X,phiY)}%
). Conversely, from (\ref{eq-d-eta(phiX,phiY)}) and (\ref{eq-phi-xi}), we
have
\[
0=d\eta \left( \varphi ^{2}X,\varphi \xi \right) =-\,d\eta \left( \varphi
X,\xi \right) ,
\]%
which along with (\ref{eq-phi-eta-xi}), when used in (\ref%
{eq-d-eta(phiX,phiY)}) yields
\[
d\eta (X,\varphi Y)=-\,d\eta (\varphi X,\varphi ^{2}Y)=-\,d\eta (\varphi
X,Y),
\]%
which in view of (\ref{eq-N2(X,Y)}) proves that $\overset{2}{N}=0$. $%
\blacksquare $ \medskip

From the definition of $\overset{3}{N}$ and $\overset{4}{{N}}$, it follows
that \cite[Theorem 3.1]{Sato-76} the tensor $\overset{3}{N}$ (resp. $\overset%
{4}{{N}}$) vanishes identically if and only if $\varphi $ (resp. $\eta $) is
invariant under the transformation generated by infinitesimal
transformations $\xi $. Consequently, in a normal almost paracontact
manifold, $\varphi $ and $\eta $ are invariant under the transformation
generated by infinitesimal transformations $\xi $. \medskip

The tangent sphere bundle over a Riemannian manifold has naturally an almost
paracontact structure in which $\overset{3}{N}=0$ and $\overset{4}{{N}}=0$
\cite{Sato-76-Yamagata}. Also an almost paracontact structure $(\varphi ,\xi
,\eta )$ is said to be {\em weak-normal} \cite{BHM-81} if the almost product
structures $J_{1}=\varphi +\eta \otimes \xi $ and $J_{2}=\varphi -\eta
\otimes \xi $ are integrable. Then an almost paracontact structure is normal
if and only if it is weak normal and $\overset{4}{N}=0$.

\section{$(\protect\varepsilon )$-$s$-paracontact metric manifolds\label%
{sect-eps-spc}}

The fundamental $\left( 0,2\right) $ symmetric tensor of the $\left(
\varepsilon \right) $-almost paracontact metric structure is defined by
\begin{equation}
\Phi \left( X,Y\right) \equiv g\left( X,\varphi Y\right)  \label{eq-Phi-def}
\end{equation}%
for all $X,Y\in TM$. Also, we get
\begin{equation}
\left( \nabla _{X}\Phi \right) \left( Y,Z\right) =g\left( \left( \nabla
_{X}\varphi \right) Y,Z\right) =\left( \nabla _{X}\Phi \right) \left(
Z,Y\right) ,  \label{eq-Phi-1}
\end{equation}%
\begin{equation}
\!\!\!\!\!\!\!\!\!\!\!\!(\nabla _{X}\Phi )(\varphi Y,\varphi Z)=-\,(\nabla
_{X}\Phi )(Y,Z)+\eta (Y)(\nabla _{X}\Phi )(\xi ,Z)+\eta (Z)(\nabla _{X}\Phi
)(Y,\xi )  \label{eq-Phi-2}
\end{equation}%
for all $X,Y,Z\in TM$.

\begin{defn-new}
\label{def-eps-pcm} We say that $(\varphi ,\xi ,\eta ,g,\varepsilon )$ is an
$(\varepsilon )${\em -paracontact metric structure} if
\begin{equation}
2\, \Phi \left( X,Y\right) = \left( \nabla _{X}\eta \right) Y+\left( \nabla
_{Y}\eta \right) X,\qquad X,Y\in TM.  \label{eq-str-pcm-def}
\end{equation}%
In this case $M$ is an $\left( \varepsilon \right) ${\em -paracontact metric
manifold}.
\end{defn-new}

The condition (\ref{eq-str-pcm-def}) is equivalent to
\begin{equation}
2\,\Phi =\varepsilon {\rm \pounds }_{\xi }g,  \label{eq-str-pcm-def-1}
\end{equation}%
where, ${\rm \pounds }$ is the operator of Lie differentiation. For $%
\varepsilon =1$ and $g$ Riemannian, $M$ is the usual paracontact metric
manifold \cite{Sato-77}.

\begin{defn-new}
\label{def-eps-s-pcm} An $\left( \varepsilon \right) $-almost paracontact
metric structure $(\varphi ,\xi ,\eta ,g,\varepsilon )$ is called an $%
(\varepsilon )$-$s${\em -paracontact metric structure} if
\begin{equation}
\nabla \xi =\varepsilon \varphi .  \label{eq-s-pcm-def}
\end{equation}%
A manifold equipped with an $(\varepsilon )$-$s$-paracontact structure is
said to be $(\varepsilon )$-$s${\em -paracontact metric manifold}.
\end{defn-new}

The equation (\ref{eq-s-pcm-def}) is equivalent to
\begin{equation}
\Phi \left( X,Y\right) =g\left( \varphi X,Y\right) =\varepsilon g\left(
\nabla _{X}\xi ,Y\right) =\left( \nabla _{X}\eta \right) Y,\qquad X,Y\in TM.
\label{eq-s-pcm-def-1}
\end{equation}

We have

\begin{theorem}
An $\left( \varepsilon \right) $-almost paracontact metric manifold is an $%
(\varepsilon )$-$s$-paracontact metric manifold if and only if it is an $%
\left( \varepsilon \right) $-paracontact metric manifold such that the
structure $1$-form $\eta $ is closed.
\end{theorem}

\noindent {\bf Proof.} Let $M$ be an $(\varepsilon )$-$s$-paracontact metric
manifold. Then in view of (\ref{eq-s-pcm-def-1}) we see that $\eta $ is
closed. Consequently, $M$ is an $(\varepsilon )$-paracontact metric manifold.

Conversely, let us suppose that $M$ is an $(\varepsilon )$-paracontact
metric manifold and $\eta $ is closed. Then
\[
\Phi \left( X,Y\right) =\frac{1}{2}\left\{ \left( \nabla _{X}\eta \right)
Y+\left( \nabla _{Y}\eta \right) X\right\} =\left( \nabla _{X}\eta \right)
Y,
\]%
which implies (\ref{eq-s-pcm-def-1}). $\blacksquare $ \medskip

\begin{prop}
If in an $\left( \varepsilon \right) $-almost paracontact metric manifold
the structure $1$-form $\eta $ is closed, then
\begin{equation}
\nabla _{\xi }\xi =0.  \label{eq-cm-del-xi-xi=0}
\end{equation}
\end{prop}

\noindent {\bf Proof.} First we note that $g\left( \nabla _{X}\xi ,\xi
\right) =0$ and in particular%
\[
g\left( \nabla _{\xi }\xi ,\xi \right) =0.
\]%
If $\eta $ is closed, then for any vector $X$ orthogonal to $\xi $, we get
\[
0=2\varepsilon \,d\eta \left( \xi ,X\right) =-\,\varepsilon \eta \left(
\left[ \xi ,X\right] \right) =-\,g\left( \xi ,\left[ \xi ,X\right] \right)
=-\,g\left( \xi ,\nabla _{\xi }X\right) =g\left( \nabla _{\xi }\xi ,X\right)
,
\]%
which completes the proof. $\blacksquare $ \medskip

Using techniques similar to those introduced in \cite[Section 4]{Sasaki-80},
we give the following

\begin{ex-new}
\label{ex-eps-s-pc} Let us assume the following:%
\begin{eqnarray*}
&&a,b,c,d\in \left\{ 1,\ldots ,p\right\} ,\qquad \quad \lambda ,\mu
,\upsilon \in \left\{ 1,\ldots ,q\right\} , \\
&&i,j,k\in \left\{ 1,\ldots ,p+q\right\} ,\qquad \lambda ^{\prime
}=p+\lambda ,\qquad n=p+q+1.
\end{eqnarray*}%
Let $\theta :{\Bbb R}^{p}\times {\Bbb R}^{q}\rightarrow {\Bbb R}$\ be a
smooth function. Define a\ function $\psi :{\Bbb R}^{n}\rightarrow {\Bbb R}$
by
\[
\psi \left( x^{1},\ldots ,x^{n}\right) \equiv \theta \left( x^{1},\ldots
,x^{p+q}\right) +x^{n}.
\]%
Now, define a\ $1$-form $\eta $ on ${\Bbb R}^{n}$\ by
\begin{equation}
\eta _{i}=\frac{\partial \theta }{\partial x^{i}}\equiv \theta _{i}\,,\qquad
\eta _{n}=1.  \label{eq-eta-Bucki-1}
\end{equation}%
Next, define a\ vector field $\xi $ on ${\Bbb R}^{n}$\ by
\begin{equation}
\xi \equiv \frac{\partial \quad }{\partial x^{n}}  \label{eq-xi-Bucki}
\end{equation}%
and a $\left( 1,1\right) $\ tensor field $\varphi $\ on ${\Bbb R}^{n}$\ by
\begin{equation}
\varphi X\equiv X^{a}\frac{\partial \quad }{\partial x^{a}}-X^{\lambda
^{\prime }}\frac{\partial \quad }{\partial x^{\lambda ^{\prime }}}+\left(
-\,\theta _{a}X^{a}+\theta _{\lambda ^{\prime }}X^{\lambda ^{\prime
}}\right) \frac{\partial \quad }{\partial x^{n}}.  \label{eq-phi-Bucki}
\end{equation}%
for all vector fields
\[
X=X^{a}\frac{\partial \quad }{\partial x^{a}}+X^{\lambda ^{\prime }}\frac{%
\partial \quad }{\partial x^{\lambda ^{\prime }}}+X^{n}\frac{\partial \quad
}{\partial x^{n}}.
\]%
Let $f_{i}:{\Bbb R}^{n}\rightarrow {\Bbb R}$ be $\left( p+q\right) $ smooth
functions. We define a tensor field $g$ of type $\left( 0,2\right) $ by
\[
g\left( X,Y\right) \equiv \left( f_{i}-\left( \theta _{i}\right) ^{2}\right)
X^{i}Y^{i}-\theta _{i}\theta _{j}X^{i}Y^{j}-\theta _{i}\left(
X^{i}Y^{n}+X^{n}Y^{i}\right) -X^{n}Y^{n},
\]%
where $f_{i}:{\Bbb R}^{n}\rightarrow {\Bbb R}$ are $\left( p+q\right) $
smooth functions such that
\[
f_{i}-\left( \theta _{i}\right) ^{2}>0,\qquad i\in \left\{ 1,\ldots
,p+q\right\} .
\]%
Then $\left( \varphi ,\xi ,\eta ,g\right) $ is a timelike Lorentzian almost
paracontact structure on ${\Bbb R}^{n}$. Moreover, if the $\left( p+q\right)
$ smooth functions $f_{i}:{\Bbb R}^{n}\rightarrow {\Bbb R}$ are given by
\[
f_{a}=F_{a}\left( x^{1},\ldots ,x^{p+q}\right) e^{-2x^{n}}+\left( \theta
_{a}\right) ^{2},\qquad a\in \left\{ 1,\ldots ,p\right\} ,
\]%
\[
f_{\lambda ^{\prime }}=F_{\lambda ^{\prime }}\left( x^{1},\ldots
,x^{p+q}\right) e^{2x^{n}}+\left( \theta _{\lambda ^{\prime }}\right)
^{2},\qquad \lambda \in \left\{ 1,\ldots ,q\right\} ,
\]%
for some smooth functions $F_{i}>0$, then we get a timelike Lorentzian $s$%
-paracontact manifold.
\end{ex-new}

\section{$\left( \protect\varepsilon \right) $-para Sasakian manifolds \label%
{sect-eps-PS}}

We begin with the following:

\begin{defn-new}
\label{def-eps-pS} An $\left( \varepsilon \right) $-almost contact metric
structure is called an $\left( \varepsilon \right) ${\em -para Sasakian
structure} if
\begin{equation}
\left( \nabla _{X}\varphi \right) Y=-\, g(\varphi X,\varphi Y)\xi
-\varepsilon \eta \left( Y\right) \varphi ^{2}X,\qquad X,Y\in TM,
\label{eq-eps-PS-def-1}
\end{equation}%
where $\nabla $ is the Levi-Civita connection with respect to $g$. A
manifold endowed with an $\left( \varepsilon \right) $-para Sasakian
structure is called an $\left( \varepsilon \right) ${\em -para Sasakian
manifold}.
\end{defn-new}


For $\varepsilon =1$ and $g$ Riemannian, $M$ is the usual para Sasakian
manifold \cite{Sato-77,Sasaki-80}. For $\varepsilon =-1$, $g$ Lorentzian and
$\xi $ replaced by $-\xi $, $M$ becomes a Lorentzian para Sasakian manifold \cite%
{Mat-89}.

\begin{ex-new}
\label{ex-eps-PS-3d} Let ${\Bbb R}^{3}$\ be the $3$-dimensional real number
space with a coordinate system $\left( x,y,z\right) $. We define
\[
\eta =dz\ ,\qquad \xi =\frac{\partial }{\partial z}\ ,
\]%
\[
\varphi \left( \frac{\partial }{\partial x}\right) =\frac{\partial }{%
\partial x}\ ,\qquad \varphi \left( \frac{\partial }{\partial y}\right) =-\,
\frac{\partial }{\partial y}\ ,\qquad \varphi \left( \frac{\partial }{%
\partial z}\right) =0\ ,
\]%
\[
g=e^{2\varepsilon x^{3}}\left( dx\right) ^{2}+e^{-2\varepsilon x^{3}}\left(
dy\right) ^{2}+\varepsilon \left( dz\right) ^{2}.
\]%
Then $\left( \varphi ,\xi ,\eta ,g,\varepsilon \right) $ is an $\left(
\varepsilon \right) $-para Sasakian structure.
\end{ex-new}

\begin{theorem}
An $\left( \varepsilon \right) $-para Sasakian structure $(\varphi ,\xi
,\eta ,g,\varepsilon )$ is always an $\left( \varepsilon \right) $-$s$%
-paracontact metric structure, and hence an $\left( \varepsilon \right) $%
-paracontact metric structure.
\end{theorem}

\noindent {\bf Proof.} Let $M$ be an $(\varepsilon )$-para Sasakian
manifold. Then from (\ref{eq-eps-PS-def-1}) we get
\[
\varphi \nabla _{X}\xi =-\left( \nabla _{X}\varphi \right) \xi =\varepsilon
\varphi ^{2}X,\qquad X,Y\in TM.
\]%
Operating by $\varphi $ to the above equation we get (\ref{eq-s-pcm-def}). $%
\blacksquare $ \medskip

The converse of the above Theorem is not true. Indeed, the $(\varepsilon )$-$%
s$-paracontact structure in the Example~\ref{ex-eps-s-pc} need not be $%
(\varepsilon )$-para Sasakian. \medskip


\begin{theorem}
An $\left( \varepsilon \right) $-para Sasakian structure is always normal.
\end{theorem}

\noindent {\bf Proof.} In an almost paracontact manifold $M$, we have
\begin{eqnarray}
{\overset{1}{N}}\left( X,Y\right) &=&\left( \nabla _{X}\varphi \right)
\varphi Y-\left( \nabla _{Y}\varphi \right) \varphi X+\left( \nabla
_{\varphi X}\varphi \right) Y-\left( \nabla _{\varphi Y}\varphi \right) X
\nonumber \\
&&-\, \eta \left( X\right) \nabla _{Y}\xi +\eta \left( Y\right) \nabla
_{X}\xi  \label{eq-torsion}
\end{eqnarray}%
for all vector fields $X,Y$ in $M$. Now, let $M$ be an $(\varepsilon )$-para
Sasakian manifold. Then it is $(\varepsilon )$-$s$-paracontact and therefore
using (\ref{eq-eps-PS-def-1}) and (\ref{eq-s-pcm-def}) in (\ref{eq-torsion}%
), we get ${\overset{1}{N}}=0$. $\blacksquare $ \medskip

\begin{pbm-new}
Whether a normal $\left( \varepsilon \right) $-paracontact structure is $%
(\varepsilon)$-para Sasakian or not.
\end{pbm-new}

\begin{lem}
\label{lem-eps-PS-R-1} Let $M$ be an $\left( \varepsilon \right) $-para
Sasakian manifold. Then the curvature tensor $R$ satisfies
\begin{equation}
R\left( X,Y\right) \xi =\eta \left( X\right) Y-\eta \left( Y\right) X,\qquad
X,Y\in TM.  \label{eq-eps-PS-R(X,Y)xi}
\end{equation}%
Consequently,
\begin{equation}
R\left( X,Y,Z,\xi \right) =-\,\eta \left( X\right) g\left( Y,Z\right) +\eta
\left( Y\right) g\left( X,Z\right) ,  \label{eq-eps-PS-R(X,Y,Z,xi)}
\end{equation}%
\begin{equation}
\eta \left( R\left( X,Y\right) Z\right) =-\,\varepsilon \eta \left( X\right)
g\left( Y,Z\right) +\varepsilon \eta \left( Y\right) g\left( X,Z\right) ,
\label{eq-eps-PS-eta(R(X,Y),Z)}
\end{equation}%
\begin{equation}
R\left( \xi ,X\right) Y=-\,\varepsilon g\left( X,Y\right) \xi +\eta \left(
Y\right) X  \label{eq-eps-PS-R(xi,X)Y}
\end{equation}%
for all vector fields $X,Y,Z$ in $M$.
\end{lem}

\noindent {\bf Proof.} Using (\ref{eq-s-pcm-def}), (\ref{eq-eps-PS-def-1})
and (\ref{eq-phi-eta-xi}) in
\[
R\left( X,Y\right) \xi =\nabla _{X}\nabla _{Y}\xi -\nabla _{Y}\nabla _{X}\xi
-\nabla _{\left[ X,Y\right] }\xi
\]%
we obtain (\ref{eq-eps-PS-R(X,Y)xi}). $\blacksquare $ \medskip

If we put
\[
R_{0}(X,Y)W=g(Y,W)X-g(X,W)Y,\qquad X,Y,W\in TM,
\]%
then in an $\left( \varepsilon \right) $-para Sasakian manifold $M$, the
equations (\ref{eq-eps-PS-R(X,Y)xi}) and (\ref{eq-eps-PS-R(xi,X)Y}) can be
rewritten as
\begin{equation}
R(X,Y)\xi =-\,\varepsilon R_{0}(X,Y)\xi ,  \label{eq-eps-PS-R(X,Y)xi-1}
\end{equation}%
\begin{equation}
R(\xi ,X)=-\,\varepsilon R_{0}(\xi ,X),  \label{eq-eps-PS-R(xi,X)Y-1}
\end{equation}%
respectively.

\begin{lem}
\label{lem-eps-PS-R-2} In an $\left( \varepsilon \right) $-para Sasakian
manifold $M$ the curvature tensor satisfies%
\begin{eqnarray}
R\left( X,Y,\varphi Z,W\right) \!\!\! &-&\!\!\!R\left( X,Y,Z,\varphi W\right)
\nonumber \\
&=&\varepsilon \,\Phi \left( Y,Z\right) g\left( \varphi X,\varphi W\right)
-\varepsilon \,\Phi \left( X,Z\right) g\left( \varphi Y,\varphi W\right)
\nonumber \\
&&+\,\varepsilon \,\Phi \left( Y,W\right) g\left( \varphi X,\varphi Z\right)
-\varepsilon \,\Phi \left( X,W\right) g\left( \varphi Y,\varphi Z\right)
\nonumber \\
&&+\,\eta \left( Y\right) \eta \left( Z\right) g\left( X,\varphi W\right)
-\eta \left( X\right) \eta \left( Z\right) g\left( Y,\varphi W\right)
\nonumber \\
&&+\,\eta \left( Y\right) \eta \left( W\right) g\left( X,\varphi Z\right)
-\eta \left( X\right) \eta \left( W\right) g\left( Y,\varphi Z\right) ,
\label{eq-eps-PS-R-1}
\end{eqnarray}%
\begin{eqnarray}
R\left( X,Y,\varphi Z,\varphi W\right) \!\!\! &-&\!\!\!R\left( X,Y,Z,W\right)
\nonumber \\
&=&\varepsilon \,\Phi \left( Y,Z\right) \Phi \left( X,W\right) -\varepsilon
\,\Phi \left( X,Z\right) \Phi \left( Y,W\right)  \nonumber \\
&&+\,\varepsilon g\left( \varphi X,\varphi Z\right) g\left( \varphi
Y,\varphi W\right) -\varepsilon g\left( \varphi Y,\varphi Z\right) g\left(
\varphi X,\varphi W\right)  \nonumber \\
&&+\,\eta \left( Z\right) \left\{ \eta \left( Y\right) g\left( X,W\right)
-\eta \left( X\right) g\left( Y,W\right) \right\}  \nonumber \\
&&-\,\eta \left( W\right) \left\{ \eta \left( Y\right) g\left( X,Z\right)
-\eta \left( X\right) g\left( Y,Z\right) \right\} ,  \label{eq-eps-PS-R-2}
\end{eqnarray}%
\begin{equation}
R\left( X,Y,\varphi Z,\varphi W\right) =R\left( \varphi X,\varphi
Y,Z,W\right) ,  \label{eq-eps-PS-R-3}
\end{equation}%
\begin{eqnarray}
R\left( \varphi X,\varphi Y,\varphi Z,\varphi W\right) &=&R\left(
X,Y,Z,W\right)  \nonumber \\
&&+\,\eta \left( Z\right) \left\{ \eta \left( Y\right) g\left( X,W\right)
-\eta \left( X\right) g\left( Y,W\right) \right\}  \nonumber \\
&&-\,\eta \left( W\right) \left\{ \eta \left( Y\right) g\left( X,Z\right)
-\eta \left( X\right) g\left( Y,Z\right) \right\} ,  \label{eq-eps-PS-R-4}
\end{eqnarray}%
for all vector fields $X,Y,Z,W$ in $M$.
\end{lem}

\noindent {\bf Proof.} Writing the equation (\ref{eq-eps-PS-def-1})
equivalently as
\[
(\nabla _{Y}\Phi )(Z,W)=-\,\varepsilon \eta \left( Z\right) g\left( \varphi
Y,\varphi W\right) -\varepsilon \eta (W)g\left( \varphi Y,\varphi Z\right)
,\qquad Y,Z,W\in TM,
\]%
and differentiating covariantly with respect to $X$ we get
\begin{eqnarray}
&&-\,\varepsilon \left( \nabla _{X}\nabla _{Y}\Phi \right) \left( Z,W\right)
\nonumber \\
&=&\Phi \left( X,Z\right) g\left( \varphi Y,\varphi W\right) +\eta \left(
Z\right) \left( \nabla _{X}\Phi \right) \left( Y,\varphi W\right)  \nonumber
\\
&&+\,\eta \left( Z\right) g\left( \varphi \left( \nabla _{X}Y\right)
,\varphi W\right) +\eta \left( Z\right) \left( \nabla _{X}\Phi \right)
\left( \varphi Y,W\right)  \nonumber \\
&&+\,\Phi \left( X,W\right) g\left( \varphi Y,\varphi Z\right) +\eta \left(
W\right) \left( \nabla _{X}\Phi \right) \left( Y,\varphi Z\right)  \nonumber
\\
&&+\,\eta \left( W\right) g\left( \varphi \left( \nabla _{X}Y\right)
,\varphi Z\right) +\eta \left( W\right) \left( \nabla _{X}\Phi \right)
\left( \varphi Y,Z\right)  \label{eq-eps-PS-nab-X-nab-Y-Phi}
\end{eqnarray}%
for all $X,Y,Z,W\in TM$. Now using (\ref{eq-eps-PS-nab-X-nab-Y-Phi}) in the
Ricci identity
\[
\left( \left( \nabla _{X}\nabla _{Y}-\nabla _{Y}\nabla _{X}-\nabla _{\left[
X,Y\right] }\right) \Phi \right) \left( Z,W\right) =-\,\Phi \left( R\left(
X,Y\right) Z,W\right) -\Phi \left( Z,R\left( X,Y\right) W\right)
\]%
we obtain (\ref{eq-eps-PS-R-1}). The equation (\ref{eq-eps-PS-R-2}) follows
from (\ref{eq-eps-PS-R-1}) and (\ref{eq-eps-PS-R(X,Y,Z,xi)}). The equation (%
\ref{eq-eps-PS-R-3}) follows from (\ref{eq-eps-PS-R-2}). Finally, the
equation (\ref{eq-eps-PS-R-4}) follows from (\ref{eq-eps-PS-R-2}) and (\ref%
{eq-eps-PS-R-3}). $\blacksquare $ \medskip

The equation (\ref{eq-eps-PS-R(X,Y)xi}) may also be obtained by (\ref%
{eq-eps-PS-R-4}). The equations (\ref{eq-eps-PS-R-1})-(\ref{eq-eps-PS-R-4})
are generalizations of the equations (3.2) and (3.3) in \cite%
{Mishra-82-Progr}. Now, we prove the following: \medskip

\begin{theorem}
\label{th-eps-PS-flat} An $\left( \varepsilon \right) $-para Sasakian
manifold can not be flat.
\end{theorem}

\noindent {\bf Proof.} Let $M$ be a flat $\left( \varepsilon \right) $-para
Sasakian manifold. Then from (\ref{eq-eps-PS-R(X,Y,Z,xi)}) we get
\[
\eta \left( X\right) g\left( Y,Z\right) =\eta \left( Y\right) g\left(
X,Z\right) ,
\]%
from which we obtain
\[
g\left( \varphi X,\varphi Z\right) =0
\]%
for all $X,Z\in TM$, a contradiction. $\blacksquare $ \medskip

A non-flat semi-Riemannian manifold $M$ is said to be {\em recurrent} \cite%
{Ruse-49} if its Ricci tensor $R$ satisfies the recurrence condition
\begin{equation}
\left( \nabla _{W}R\right) \left( X,Y,Z,V\right) =\alpha \left( W\right)
R\left( X,Y,Z,V\right) ,\quad X,Y,Z,V\in TM,  \label{eq-rec-sp}
\end{equation}%
where $\alpha $ is a $1$-form. If $\alpha =0$ in the above equation, then
the manifold becomes {\em symmetric} in the sense of Cartan \cite{Cartan-26}%
. We say that $M$ is proper recurrent, if $\alpha \neq 0$.

\begin{theorem}
\label{th-eps-PS-recurrent} An $\left( \varepsilon \right) $-para Sasakian
manifold cannot be proper recurrent.
\end{theorem}

\noindent {\bf Proof.} Let $M$ be a recurrent $\left( \varepsilon \right) $%
-para Sasakian manifold. Then from (\ref{eq-rec-sp}), (\ref%
{eq-eps-PS-R(X,Y,Z,xi)}) and (\ref{eq-s-pcm-def}) we obtain
\begin{eqnarray}
\varepsilon R\left( X,Y,Z,\varphi W\right) &=&g\left( X,Z\right) \left\{
\Phi \left( Y,W\right) -\alpha \left( W\right) \eta \left( Y\right) \right\}
\nonumber \\
&&-g\left( Y,Z\right) \left\{ \Phi \left( X,W\right) -\alpha \left( W\right)
\eta \left( X\right) \right\}  \label{eq-eps-PS-rec}
\end{eqnarray}%
for all $X,Y,Z,W\in TM$. Putting $Y=\xi $ in the above equation, we get
\[
\alpha \left( W\right) g\left( \varphi X,\varphi Z\right) =0,\qquad X,Z,W\in
TM,
\]%
a contradiction. $\blacksquare $ \medskip

Let $n\geq 2$ and $0\leq \nu \leq n$. Then \cite[Definition 23, p. 110]%
{ONeill-83}:

\begin{enumerate}
\item The {\em pseudosphere} of radius $r>0$ in $R_{\nu }^{n+1}$ is the
hyperquadric
\[
S_{\nu }^{n}\left( r\right) =\left\{ p\in R_{\nu }^{n+1}:\left\langle
p,p\right\rangle =r^{2}\right\}
\]%
with dimension $n$ and index $\nu $.

\item The {\em pseudohyperbolic space} of radius $r>0$ in $R_{\nu +1}^{n+1}$
is the hyperquadric
\[
H_{\nu }^{n}\left( r\right) =\left\{ p\in R_{\nu +1}^{n+1}:\left\langle
p,p\right\rangle =-\,r^{2}\right\}
\]%
with dimension $n$ and index $\nu $.
\end{enumerate}

\begin{theorem}
\label{th-eps-PS-sym} An $\left( \varepsilon \right) $-para Sasakian
manifold is symmetric if and only if it is of constant curvature $%
-\,\varepsilon $. Consequently, a symmetric spacelike $\left( \varepsilon
\right) $-para Sasakian manifold is locally isometric to a pseudohyperbolic
space $H_{\nu }^{n}(1)$ and a symmetric timelike $\left( \varepsilon \right)
$-para Sasakian manifold is locally isometric to a pseudosphere $S_{\nu
}^{n}(1)$.
\end{theorem}

\noindent {\bf Proof.} Let $M$ be a symmetric $\left( \varepsilon \right) $%
-para Sasakian manifold. Then putting $\alpha =0$ in (\ref{eq-eps-PS-rec})
we obtain
\[
\varepsilon R\left( X,Y,Z,\varphi W\right) =g\left( X,Z\right) \Phi \left(
Y,W\right) -g\left( Y,Z\right) \Phi \left( X,W\right)
\]%
for all $X,Y,Z,W\in TM$. Writing $\varphi W$ in place of $W$ in the above
equation and using (\ref{eq-metric-1}) and (\ref{eq-eps-PS-R(X,Y,Z,xi)}), we
get
\begin{equation}
R\left( X,Y,Z,W\right) =-\,\varepsilon \left\{ g\left( Y,Z\right) g\left(
X,W\right) -g\left( X,Z\right) g\left( Y,W\right) \right\} ,
\label{eq-eps-PS-sp-form}
\end{equation}%
which shows that $M$ is a space of constant curvature $-\,\varepsilon $. The
converse is trivial. $\blacksquare $ \medskip

\begin{cor}
\label{cor-eps-PS-const-curv-2} If an $\left( \varepsilon \right) $-para
Sasakian manifold is of constant curvature, then
\begin{eqnarray}
&&\Phi \left( Y,Z\right) \Phi \left( X,W\right) -\Phi \left( X,Z\right) \Phi
\left( Y,W\right)  \nonumber \\
&=&-\,g\left( \varphi Y,\varphi Z\right) g\left( \varphi X,\varphi W\right)
+g\left( \varphi X,\varphi Z\right) g\left( \varphi Y,\varphi W\right)
\label{eq-eps-PS-const-curv-2}
\end{eqnarray}%
for all $X,Y,Z,W\in TM$.
\end{cor}

\noindent {\bf Proof.} Obviously, if an $\left( \varepsilon \right) $-para
Sasakian manifold is of constant curvature $k$, then $k=-\,\varepsilon $.
Therefore, using (\ref{eq-eps-PS-sp-form}) in (\ref{eq-eps-PS-R-2}) we get (%
\ref{eq-eps-PS-const-curv-2}). $\blacksquare $ \medskip

Apart from recurrent spaces, semi-symmetric spaces are another well-known
and important natural generalization of symmetric spaces. A semi-Riemannian
manifold $\left( M,g\right) $ is a {\em semi-symmetric space} if its
curvature tensor $R$ satisfies the condition
\[
R(X,Y)\cdot R=0
\]%
for all vector fields $X,Y$ on $M$, where $R\left( X,Y\right) $ acts as a
derivation on $R$. Symmetric spaces are obviously semi-symmetric, but the
converse need not be true. In fact, in dimension greater than two there
always exist examples of semi-symmetric spaces which are not symmetric. For
more details we refer to \cite{BKV-96}. \medskip

Given a class of semi-Riemannian manifolds, it is always interesting to know
that whether, inside that class, semi-symmetry implies symmetry or not.
Here, we prove the following:

\begin{theorem}
\label{th-eps-PS-semi-sym} In an $\left( \varepsilon \right) $-para Sasakian
manifold, the condition of semi-symmetry implies the condition of symmetry.
\end{theorem}

\noindent {\bf Proof.} Let $M$ be a symmetric $\left( \varepsilon \right) $%
-para Sasakian manifold. Let the condition of being semi-symmetric be true,
that is,
\[
R\left( V,U\right) \cdot R=0,\qquad V,U\in TM.
\]%
In particular, from the condition $R(\xi ,U)\cdot R=0$, we get
\[
0=\left[ R(\xi ,U),R(X,Y)\right] \xi -R(R(\xi ,U)X,Y)\xi -R(X,R(\xi ,U)Y)\xi
,
\]%
which in view of (\ref{eq-eps-PS-R(xi,X)Y-1}) gives
\begin{eqnarray*}
0 &=&g(U,R(X,Y)\xi )\xi -\eta (R(X,Y)\xi )U \\
&&-\ g(U,X)R(\xi ,Y)\xi +\eta (X)R(U,Y)\xi -g(U,Y)R(X,\xi )\xi \\
&&+\ \eta (Y)R(X,U)\xi -\eta (U)R(X,Y)\xi +R(X,Y)U.
\end{eqnarray*}%
Equation (\ref{eq-eps-PS-R(X,Y)xi-1}) then gives
\[
R=-\,\varepsilon R_{0}.
\]%
Therefore $M$ is of constant curvature $-\,\varepsilon $, and hence
symmetric. $\blacksquare $ \medskip

In view of Theorem~\ref{th-eps-PS-sym} and Theorem~\ref{th-eps-PS-semi-sym},
we have the following:

\begin{cor}
Let $M$ be an $\left( \varepsilon \right) $-para Sasakian manifold. Then the
following statements are equivalent:

\begin{enumerate}
\item[{\bf (a)}] $M$ is symmetric.

\item[{\bf (b)}] $M$ is of constant curvature $-\,\varepsilon $.

\item[{\bf (c)}] $M$ is semi-symmetric.

\item[{\bf (d)}] $M$ satisfies $R\left( \xi ,U\right) \cdot R=0$.
\end{enumerate}
\end{cor}

Now, we need the following:

\begin{lem}
\label{lem-eps-PS-Ricci} In an $n$-dimensional $\left( \varepsilon \right) $%
-para Sasakian manifold $M$ the Ricci tensor $S$ satisfies
\begin{equation}
S\left( \varphi Y,\varphi Z\right) =S\left( Y,Z\right) +\left( n-1\right)
\eta \left( Y\right) \eta \left( Z\right)  \label{eq-eps-PS-Ricci-1}
\end{equation}%
for all $Y,Z\in TM$. Consequently,
\begin{equation}
S\left( \varphi Y,Z\right) =S\left( Y,\varphi Z\right) ,
\label{eq-eps-PS-Ricci-2}
\end{equation}%
\begin{equation}
S\left( Y,\xi \right) =-\left( n-1\right) \eta \left( Y\right) .
\label{eq-eps-PS-S(Y,xi)}
\end{equation}
\end{lem}

\noindent {\bf Proof.} Contracting the equation (\ref{eq-eps-PS-R-4}), we
get (\ref{eq-eps-PS-Ricci-1}). Replacing $Z$ by $\varphi Z$ in (\ref%
{eq-eps-PS-Ricci-1}) we get (\ref{eq-eps-PS-Ricci-2}). Putting $Z=\xi $ in (%
\ref{eq-eps-PS-Ricci-1}) we get (\ref{eq-eps-PS-S(Y,xi)}). $\blacksquare $
\medskip

A semi-Riemannian manifold $M$ is said to be {\em Ricci-recurrent} \cite%
{Pat-52} if its Ricci tensor $S$ satisfies the condition
\begin{equation}
\left( \nabla _{X}S\right) \left( Y,Z\right) =\alpha \left( X\right) S\left(
Y,Z\right) ,\qquad X,Y,Z\in TM,  \label{eq-Ricci-rec}
\end{equation}%
where $\alpha $ is a $1$-form. If $\alpha =0$ in the above equation, then
the manifold becomes {\em Ricci-symmetric}. We say that $M$ is proper
Ricci-recurrent, if $\alpha \neq 0$.

\begin{theorem}
\label{th-eps-PS-Ricci-rec} An $\left( \varepsilon \right) $-para Sasakian
manifold can not be proper Ricci-recurrent.
\end{theorem}

\noindent {\bf Proof.} Let $M$ be an $n$-dimensional $\left( \varepsilon
\right) $-para Sasakian manifold. If possible, let $M$ be proper
Ricci-recurrent. Then
\begin{equation}
\left( \nabla _{X}S\right) \left( Y,\xi \right) =\alpha \left( X\right)
S\left( Y,\xi \right) =-\,\left( n-1\right) \alpha \left( X\right) \eta
\left( Y\right) .  \label{eq-eps-PS-Ricci-rec-1}
\end{equation}%
But we have
\begin{equation}
\left( \nabla _{X}S\right) \left( Y,\xi \right) =\left( n-1\right) \left(
\nabla _{X}\eta \right) Y-\varepsilon S\left( Y,\varphi X\right) .
\label{eq-eps-PS-Ricci-rec-2}
\end{equation}%
Using (\ref{eq-eps-PS-Ricci-rec-2}) in (\ref{eq-eps-PS-Ricci-rec-1}) we get
\begin{equation}
\varepsilon S\left( \varphi X,Y\right) +\left( n-1\right) \Phi \left(
X,Y\right) =\left( n-1\right) \alpha \left( X\right) \eta \left( Y\right)
\label{eq-eps-PS-Ricci-rec-3}
\end{equation}%
Putting $Y=\xi $ in the above equation, we get $\alpha \left( X\right) =0$,
a contradiction. $\blacksquare $ \medskip

A semi-Riemannian manifold $M$ is said to be {\em Ricci-semi-symmetric} \cite%
{Mir-92} if its Ricci tensor $S$ satisfies the condition
\[
R(X,Y)\cdot S=0
\]%
for all vector fields $X,Y$ on $M$, where $R\left( X,Y\right) $ acts as a
derivation on $S$. \medskip

In last, we prove the following:

\begin{theorem}
\label{th-eps-PS-Ricci} For an $n$-dimensional $\left( \varepsilon \right) $%
-para Sasakian manifold $M$, the following three statements are equivalent:

\begin{enumerate}
\item[{\bf (a)}] $M$ is an Einstein manifold.

\item[{\bf (b)}] $M$ is Ricci-symmetric.

\item[{\bf (c)}] $M$ is Ricci-semi-symmetric.
\end{enumerate}
\end{theorem}

\noindent {\bf Proof.} Obviously, the statement {\bf (a)} implies each of
the statements {\bf (b)} and {\bf (c)}. Let {\bf (b)} be true. Then putting $%
\alpha =0$ in (\ref{eq-eps-PS-Ricci-rec-3}) we get
\begin{equation}
\varepsilon S\left( \varphi X,Y\right) +\left( n-1\right) \Phi \left(
X,Y\right) =0.  \label{eq-eps-PS-Ricci-sym-1}
\end{equation}%
Replacing $X$ by $\varphi X$ in the above equation, we get
\begin{equation}
S=-\,\varepsilon \left( n-1\right) g,  \label{eq-eps-PS-Ricci-sym-2}
\end{equation}%
which shows that the statement {\bf (a)} is true. In last, let {\bf (c)} be
true. In particular,
\[
\left( R(\xi ,X)\cdot S\right) \left( Y,\xi \right) =0
\]%
implies that
\[
S\left( R(\xi ,X)Y,\xi \right) +S\left( Y,R(\xi ,X)\xi \right) =0,
\]%
which in view of (\ref{eq-eps-PS-R(xi,X)Y}) and (\ref{eq-eps-PS-S(Y,xi)})
again gives (\ref{eq-eps-PS-Ricci-sym-2}). This completes the proof. $%
\blacksquare $ \medskip

\noindent {\bf Acknowledgement.} This paper was prepared during the first
visit of the first author to \.{I}n\"{o}n\"{u} University, Turkey in July
2009. The first author was supported by the Scientific and Technical
Research Council of Turkey (T\"{U}B\.{I}TAK) for Advanced Fellowships
Programme.

\noindent Mukut Mani Tripathi

\noindent Department of Mathematics, Banaras Hindu University

\noindent Varanasi 221 005, India.

\noindent Email: mmtripathi66@@yahoo.com

\smallskip

\noindent Erol K\i l\i \c{c}

\noindent Department of Mathematics, Faculty of Arts and Sciences, \.{I}n%
\"{o}n\"{u} University

\noindent 44280 Malatya, Turkey

\noindent Email: ekilic@@inonu.edu.tr

\smallskip

\noindent Selcen Y\"{u}ksel Perkta\c{s}

\noindent Department of Mathematics, Faculty of Arts and Sciences, \.{I}n%
\"{o}n\"{u} University

\noindent 44280 Malatya, Turkey

\noindent Email: selcenyuksel@@inonu.edu.tr

\smallskip

\noindent Sad\i k Kele\c{s}

\noindent Department of Mathematics, Faculty of Arts and Sciences, \.{I}n%
\"{o}n\"{u} University

\noindent 44280 Malatya, Turkey

\noindent Email: keles@@inonu.edu.tr

\end{document}